\documentclass[final,10.5pt]{amsart}
\usepackage{amsmath}
\usepackage{amssymb}
\usepackage{amscd}
\usepackage{xspace}
\usepackage{verbatim}
\usepackage[notcite,notref]{showkeys}

\newcommand{\fnz}{\footnotesize}

\newcommand{\be}{\begin{equation}}
\newcommand{\bel}[1]{\begin{equation}\label{#1}}
\newcommand{\ee}{\end{equation}}
\newtheorem{subn}{\name}

\newcommand{\bsn}[1]{\def\name{#1}\begin{subn}}
\newcommand{\esn}{\end{subn}}
\newtheorem{sub}{\name}
\newcommand{\bs}{\begin{sub}}
\newcommand{\es}{\end{sub}}
\newcommand{\bsl}[1]{\begin{sub}\label{#1}}
\newcommand{\bth}[1]{\def\name{Theorem}\begin{sub}\label{t:#1}}
\newcommand{\blemma}[1]{\def\name{Lemma}\begin{sub}\label{l:#1}}
\newcommand{\bcor}[1]{\def\name{Corollary}\begin{sub}\label{c:#1}}
\newcommand{\bdef}[1]{\def\name{Definition}\begin{sub}\label{d:#1}}
\newcommand{\bprop}[1]{\def\name{Proposition}\begin{sub}\label{p:#1}}

\newcommand{\rth}[1]{Theorem~\ref{t:#1}}

\newcommand{\BA}{\begin{array}}
\newcommand{\EA}{\end{array}}
\newcommand{\BAN}{\renewcommand{\arraystretch}{1.2}
\setlength{\arraycolsep}{2pt}\begin{array}}
\newcommand{\BAV}[2]{\renewcommand{\arraystretch}{#1}
\setlength{\arraycolsep}{#2}\begin{array}}
\newcommand{\BSA}{\begin{subarray}}
\newcommand{\ESA}{\end{subarray}}
\newcommand{\BAL}{\begin{aligned}}
\newcommand{\EAL}{\end{aligned}}
\newcommand{\BALG}{\begin{alignat}}
\newcommand{\EALG}{\end{alignat}}
\newcommand{\BALGN}{\begin{alignat*}}
\newcommand{\EALGN}{\end{alignat*}}
\newcommand{\note}[1]{\textit{#1.}\hspace{2mm}}
\newcommand{\Proof}{\noindent\note{Proof}}

\newcommand{\forevery}{\quad \forall}

\newcommand{\abs}[1]{\left |#1\right |}
\newcommand{\norm}[1]{\left \|#1\right \|}

\newcommand{\myfrac}[2]{{\displaystyle \frac{#1}{#2} }}
\newcommand{\myint}[2]{{\displaystyle \int_{#1}^{#2}}}

\newcommand{\prt}{\partial}

\newcommand{\ti}{\times}

\newcommand{\nind}{\noindent}

\def\ga{\alpha}     \def\gb{\beta}       \def\gg{\gamma}
       \def\gd{\delta}      \def\ge{\epsilon}
\def\gth{\theta}

\def\gm{\mu}                 
    \def\gr{\rho}        
\def\gs{\sigma}       \def\gt{\tau}
      \def\gw{\omega}
                
     \def\Gd{\Delta}

\def\Gw{\Omega}              

   \def\BBN {\mathbb N}    
   \def\BBR {\mathbb R}

\begin{document}%
\nind{\small {\bf Partial Differential Equations}}/{\small {\it
\'Equations aux d\'eriv\'ees partielles}}\\

\nind{\LARGE\bf Diffusion versus absorption in semilinear parabolic
problems\footnote {To appear in{ \it C. R. Acad. Sci. Paris, Ser. I.}}}\\

\nind{\bf Andrey Shishkov}\\
{\small Institute of Applied Mathematics and Mechanics of NAS of
Ukraine,
R. Luxemburg str. 74, 83114 Donetsk, Ukraine\\
 Email: {\it shishkov@iamm.ac.donetsk.ua }}. \\[-2mm]

\nind {\bf Laurent V\'eron}\\
{\small Laboratoire de Math\'ematiques et Physique Th\'eorique, CNRS
UMR 6083, Facult\'e des Sciences, 37200 Tours, France. \\
Email: {\it veronl@univ-tours.fr}}\\
\underline{\phantom{-----------------------------------------------------------------------------------------------}}

\begin{quotation}
 {\fnz{\sc \bf{Abstract}.}
 We study the limit, when $k\to\infty$, of the solutions $u=u_{k}$ of
(E)
 $\prt_{t}u-\Gd u+ h(t)u^q=0$ in $\BBR^N\ti (0,\infty)$,
$u_{k}(.,0)=k\gd_{0}$, with $q>1$, $h(t)>0$. If
$h(t)=e^{-\gw(t)/t}$ where $\gw>0$ satisfies to
$\int_{0}^1\sqrt{\gw(t)}t^{-1}dt<\infty$, the limit function
$u_{\infty}$ is a solution of (E) with a single singularity at
$(0,0)$, while if $\gw(t)\equiv 1$, $u_{\infty}$ is the maximal
solution of (E). We examine
similar questions for equations such as  $\prt_{t}u-\Gd u^m+
h(t)u^q=0$ with $m>1$ and $\prt_{t}u-\Gd u+ h(t)e^{u}=0$.
 }\end {quotation}
\indent {\large\bf Diffusion versus absorption dans des probl\`emes
paraboliques semi-lin\'eaires}
\begin{quotation}
 {\fnz{\sc \bf{R\'esum\'e}.}
%
 \setcounter{equation}{0}
Nous \'etudions la limite, quand $k\to\infty$, des solutions
$u=u_{k}$ de (E)
 $\prt_{t}u-\Gd u+ h(t)u^q=0$ dans $\BBR^N\ti (0,\infty)$,
$u_{k}(.,0)=k\gd_{0}$ avec $q>1$, $h(t)>0$. Nous
montrons que si $h(t)=e^{-\gw(t)/t}$ o\`u $\gw>0$ v\'erifie  $\int_{0}^1\sqrt{\gw(t)}t^{-1}dt<\infty$, la
fonction limite $u_{\infty}$ est une solution of (E) avec une
singularit\'e isol\'ee en $(0,0)$, alors que si $\gw(t)\equiv 1$,
$u_{\infty}$ est la solution maximale de (E). Nous examinons des questions semblables pour des
\'equations  des type suivants $\prt_{t}u-\Gd u^m+ h(t)u^q=0$ avec
$m>1$ et $\prt_{t}u-\Gd u+ h(t)e^{u}=0$.
 }
\end
{quotation}\underline{\phantom{-----------------------------------------------------------------------------------------------}}\\

\nind{\large\bf{Version fran\c caise abr\'eg\'ee}}\medskip

Soit $q>1$ et $h:\BBR_{+}\mapsto \BBR_{+}$ une fonction continue,
croissante telle que $h(t)>0$ pour $t>0$. Il est facile de v\'erifier
que toute solution positive $u$ de
\begin {equation}\label {main1}
 \prt_{t}u-\Gd u+ h(t)u^q=0\quad\mbox {dans }\BBR^N\ti]0,+\infty[
\end {equation}
satisfait \`a
\begin {equation}\label {max1}
 u(x,t) \leq
U(t):=\left((q-1)\myint{0}{t}h(s)\,ds\right)^{-1/(q-1)}\forevery
(x,t)\in \BBR^N\ti ]0,+\infty[.
 \end {equation}
 Si $h\in L^1(0,1,\,t^{Nq/2}dt)$, il est classique que pour tout
$k>0$ il existe une unique solution (dite fondamentale) $u=u_{k}$ de
(\ref {main1}) sur $\BBR^N\ti]0,+\infty[$ v\'erifiant
$u_{k}(.,0)=k\gd_{0}$. Par le principe du maximum $k\mapsto u_{k}$
est croissant et deux cas peuvent se produire:\\[-4mm]

\noindent (i) ou bien $u_{\infty}=\lim_{k\to\infty}u_{k}=U$. {\it
Explosion initiale compl\`ete}.\\[-4mm]

\noindent (ii) ou bien $u_{\infty}$ est une solution de (\ref
{main1}) singuli\`ere en $(0,0)$
v\'erifiant
 $\lim_{t\to 0}u_{\infty}(x,t)$

 $\;\;=0$ pour tout $x\neq 0$. {\it Explosion initiale ponctuelle}. \\

\noindent {\bf Theorem 1.}{ \it (I) Si $h(t)=e^{-\gs/t}$ pour un
$\gs>0$, alors $u_{\infty}=U$.\\[-4mm]

\noindent (II) Si $h(t)=e^{-\gw(t)/t}$ o\`u $\gw$ est  monotone
croissante sur $]0,+\infty[$ et v\'erifie, pour un $\ga\in [0,1[$, $\inf\{\gw(t)/t^{\ga}:0<t\leq 1\}>0 $ et  
$\int_{0}^1\sqrt {\gw(t)}\,t^{-1}\,dt<\infty$,
 alors $u_{\infty}$ a une explosion initiale ponctuelle.
}\\

 Dans le cas de l'\' equation
  \begin {equation}\label {main2}
 \prt_{t}u-\Gd u+ h(t)e^{u}=0\quad\mbox {dans }\BBR^N\ti]0,+\infty[,
\end {equation}
toute solution $u$ satisfait \`a
\begin {equation}\label {max2}
 u(x,t) \leq \tilde U(t):=-\ln\left(\int_{0}^th(s)\,ds\right)\forevery
(x,t)\in \BBR^N\ti ]0,+\infty[,
 \end {equation}
 et l'existence d'une solution fondamentale $u=u_{k}$ est assur\'ee
si $h(t)=e^{-b(t)}$ avec
$ \lim_{t\to+\infty}t^{N/2}b(t)=+\infty$.\\

\noindent {\bf Theorem 2.}{\it (I) Si $h(t)=e^{-e^{\gs/t}}$ pour un
$\gs>0$, alors $u_{\infty}=\tilde U$.\\[-4mm]

\noindent (II) Si $h(t)=e^{-e^{\gw(t)/t}}$ o\`u $\gw$ v\'erifie
les conditions du Th\'eor\`eme 1,
alors $u_{\infty}$ a une explosion initiale ponctuelle.
}\\

Nos m\'ethodes nous permettent ausi de traiter l'\'equation des
 milieux poreux avec absorption.\\ [-2mm]

 \nind{\large\bf{Main results}}\\ [-2mm]
\setcounter{equation}{0}

\noindent Let $q>1$ and $h:(0,\infty)\mapsto (0,\infty)$ be a continuous
nondecreasing function. It is easy to prove that any positive
solution $u$ of
\begin {equation}\label {main1}
 \prt_{t}u-\Gd u+ h(t)u^q=0\quad\mbox {dans }\BBR^N\ti(0,+\infty)
\end {equation}
verifies
\begin {equation}\label {max1}
 u(x,t) \leq
U(t):=\left((q-1)\myint{0}{t}h(s)\,ds\right)^{-1/(q-1)}\forevery
(x,t)\in \BBR^N\ti (0,\infty).
 \end {equation}
 If $h\in L^1(0,1,\,t^{Nq/2}dt)$, it is classical that, for any $k>0$,
there exists a unique solution (called fundamental) $u=u_{k}$ of
(\ref {main1}) sur $\BBR^N\ti(0,\infty)$ such that
$u_{k}(.,0)=k\gd_{0}$. By the maximum principle $k\mapsto u_{k}$
is
increasing and  the following alternative occurs:\\[-4mm]

\noindent (i) either $u_{\infty}=\lim_{k\to\infty}u_{k}=U$. {\it
Complete initial blow-up}.\\[-4mm]

\noindent (ii) or $u_{\infty}$ is a solution of (\ref {main1})
singular at $(0,0)$  such that

 $\;\;\lim_{t\to 0}u_{\infty}(x,t)=0$ for all $x\neq 0$. {\it
Single-point initial blow-up}.
 \bth {th1} (I) If $h(t)=e^{-\gs/t}$ for some $\gs>0$, then
$u_{\infty}=U$.\\[-4mm]

\noindent (II) If $h(t)=e^{-\gw(t)/t}$ where $\gw$ is nondecreasing
 on $(0,+\infty )$ and verifies, for some $\ga\in [0,1)$, $\inf\{\gw(t)/t^{\ga}:0<t\leq 1\}>0 $ and
\begin {equation}\label {dini}
\int_{0}^1\myfrac {\sqrt {\gw(t)}\,dt}{t}<+\infty,
\end {equation}
then $u_{\infty}$ has single-point initial blow-up.
 \es
Concerning  equation
  \begin {equation}\label {main2}
 \prt_{t}u-\Gd u+ h(t)e^{u}=0\quad\mbox {dans }\BBR^N\ti(0,+\infty),
\end {equation}
any solution $u$ verifies
\begin {equation}\label {max2}
 u(x,t) \leq \tilde U(t):=-\ln\left(\int_{0}^th(s)\,ds\right)\forevery
(x,t)\in \BBR^N\ti (0,+\infty).
 \end {equation}
 and the existence of a fundamental solution $u=u_{k}$ is ensured if
$h(t)=e^{-b(t)}$ where
$ \lim_{t\to+\infty}t^{N/2}b(t)=+\infty$.
 \bth {th2} (I) If $h(t)=O(e^{-e^{\gs/t}})$ for some $\gs>0$, then
$u_{\infty}=\tilde U$.\\[-4mm]

\noindent (II) If $h(t)=e^{-e^{\gw(t)/t}}$ where $\gw$ satisfies the conditions of \rth {th1}, then $u_{\infty}$ has single-point initial blow-up.
\es
Our methods apply to equations of porous media type
\begin{equation}\label{main3}
\prt_{t}u-\Gd u^m+ h(t)u^q=0\quad \mbox {in }\BBR^N
\ti (0,\infty),
\end{equation}
with $m>1$, $q>1$ and $h:(0,\infty)\mapsto (0,\infty)$ is
nondecreasing. As above, any positive solution satisfies
(\ref{max1}). If $h\in L^1((0,1;t^{-(q-1)/(m-1+2N^{-1})}dt)$, for any
$k>0$ there exists a solution
$u=u_{k}$ of (\ref {main3}) such that $u_{k}(.,0)=k\gd_{0}$. Since
$k\mapsto u_{k}$ is increasing, the same alternative as in case of
(\ref {main1}) occurs concerning $u_{\infty}$.
 \bth {th3} Assume $q>m>1$. (I) If $h(t)=O(t^{(q-m)/(m-1)})$ , then
$u_{\infty}=U$.\\[-4mm]

\noindent (II) If $h(t)=t^{(q-m)/(m-1)}\gw^{-1}(t)$ where $\gw$ is
nondecreasing and positive on $(0,+\infty )$ and verifies
\begin{equation}\label{dini2}
\int_{0}^1\myfrac {\gw^\gth (t)\,dt}{t}<+\infty,
\end{equation}
where
$$\gth=\myfrac {m^2-1}{\left( N(m-1)+2(m+1)\right)(q-1)},
$$
then $u_{\infty}$ has single-point initial blow-up.
 \es

 \noindent {\it Sketch of the proofs.} The complete initial blow-up
results are proved by constructing local subsolutions by modifying
the very singular solutions of  some related equations. Since for
equation (\ref {main1}),  the proof is already given in \cite {MV1} we
shall outline the (more complicated) construction for equation (\ref {main2}).
 \blemma{L1} If $h(t)=\gs t^{-2}e^{\gs t^{-1}-e^{-\gs/t}}$ for some
$\gs>0$, complete initial blow-up occurs for equation (\ref {main2}).
 \es
 \Proof Writing $h(t)=e^{-a(t)}$ is is first observed that
fundamental solutions $u_{k}$ of (\ref {main2})
exist for all $k>0$ if $\lim_{{t\to 0}}t^{N/2}a(t)=\infty$. For
$\ell>1$, let $v=v_{\infty,\ell}$ be the very singular solution of
 \begin {equation}\label {vss}
\prt_{t}v-\Gd v+ct^{\ga_{\ell}}v^\ell=0
\end {equation}
 in $\BBR^N\ti (0,\infty)$, where $\ga_{\ell}$ and $c$ are positive
constants. The choice of $\ga_{\ell}=(N+2)/(\ell -1)/2-1$ ensures the
existence of $v_{\infty,\ell}$. Furthermore, if we write
 $$v_{\infty,\ell}(x,t)=\left(\myfrac
{2c}{N+2}\right)^{1/(\ell-1)}t^{-(1+N/2)} f_{\ell}(x/\sqrt {t}),
 $$
 then  $f_{\ell}(\eta)\leq 1$ for $\eta\in \BBR^N$ and
  \begin {equation}\label {vss2}
\Gd f_{\ell}+\myfrac {1}{2}Df_{\ell}.\eta
+\myfrac{N+2}{2}f_{\ell}-f_{\ell}^\ell=0.
\end {equation}
By the maximum principle $0<f_{\ell}<f_{\ell'}\leq 1$ for
$\ell'>\ell>1$. For the particular choice $\ell^*=(N+4)/(N2)$, we can
use the expression of the asymptotic expansion of the very singular
solution given in \cite {BPT},
$$f_{\ell^*}(\eta)=C{\abs \eta}^2e^{-{\abs \eta}^2/4}(1+\circ
(1))\mbox { as }{\abs \eta}\to\infty,
$$
from which follows $f_{\ell}(\eta)\geq
f_{\ell^*}(\eta)\geq\gd^*({\abs \eta}^2+1)e^{-{\abs \eta}^2/4}$ for
some $\gd^*>0$, any $\eta\in\BBR^N$ and $\ell\geq \ell^*$. Thus there
exists $\gd>0$ depending only on $N$ such that
  \begin {equation}\label {vss3}
v_{\infty,\ell}(x,t)\geq \gd c^{1/(\ell-1)}t^{-1-N/2}({\abs
x}^2+t)e^{-{\abs x}^2/4t}
\forevery (x,t)\in\BBR^N\ti(0,\infty).
\end {equation}
Because any positive solution $u$ of (\ref {main2}) satisfies
(\ref{max2}), we have to prove that we can fix $c$ and $\gt>0$ such
that
  \begin {equation}\label {vss4}
  ct^{\ga_{\ell}}(\gr^\ell+1)\geq h(t)e^\gr\forevery (t,\gr)\in
(0,\gt]\ti [0,\tilde V(t)].
\end {equation}
Writing $h$ under the form $h(t)=-\gw'(t)e^{\gw(t)}$ where
$\gw(t)=e^{\gg(t)}$ and $\gg$ is a positive decreasing $C^1$
function, infinite at $t=0$, we first notice that it is sufficient to
prove this inequality for $\gr=\tilde U(t)$, and in that case
  \begin {equation}\label {vss5}
  ct^{\ga_{\ell}}(e^{\ell\gg(t)}+1)\geq-\gg'(t)e^{\gg(t)}\forevery
t\in (0,\gt].
\end {equation}
We take now $\gg(t)=\gs/ t$, and prove that there exists $\gb>0$,
depending only on $N$ such that,
for any $0<\gt\leq\gb\gs$, estimate (\ref {vss4})
holds with
$$ c=e^{(1-\ell)\gs/\gt-2^{-1}(\ell (N+2)-N)\ln\gt}.
$$
The maximum principle and (\ref {vss4}) imply that for any $\ell>1$
and $k>0$ the solutions $u=u_{k}$ of (\ref {main2}) and $v=\tilde
v_{k}$ of
$$
\prt_{t}v-\Gd v+ct^{\ga_{\ell}}(v^\ell+1)=0
$$
with initial data $k\gd_{0}$ verifies $0\leq \tilde v_{k,\ell}\leq
u_{k}$, on $(0,\gt]$. Therefore $v_{\infty,\ell}\leq
u_{\infty}+ct^{\ga_{\ell}+1}/(\ga_{\ell}+1)$ on $(0,\gt]$ leads to
$$
u_{\infty}(x,\gt)\geq
\gd \gt^{-1-N/2}({\abs x}^2+\gt)
e^{\frac {4\gs-{\abs x}^2}{4\gt}-\frac {\ell
(N+2)-N}{2(1-\ell)}\ln\gt}
$$
Thus $\lim_{\gt\to 0}u_{\infty}(x,\gt)=\infty$, locally uniformly in
$B_{2\sqrt\gs}$, which implies the result.
\\

The proof of \rth {th2} follows from the fact that for any
$\gs>\gs'>0$ there exists an interval $(0,\gth]$ where
$\gs' t^{-2}e^{\gs' t^{-1}-e^{-\gs'/t}}\geq e^{-e^{\gs/t}}$.\\

 The single-point initial blow-up is proved by local energy methods.
Because of their high degree of technicality we shall just give a
short sketch of them in the simplest case of \rth{th1}. For $k>0$,
let $u_{k}=u$ be the solution of the next result.
 \begin{equation}\label{energ1}\left\{\BA{l}
 \prt_{t} u-\Gd u+h(t)\abs u^{q-1}u=0\quad \mbox {in }\BBR^N\ti
(0,\infty)\\[2mm]
 u(x,0)=u_{0,k}(x)=M_{k}^{1/2}k^{-N/2}\eta_{k}(x)\forevery x\in\BBR^N,
 \EA\right.
 \end{equation}

\noindent where $\eta_{k}\in C(\BBR^N)$ is nonnegative, has
compact support in $B_{k^{-1}}$,  converges weakly to $\gd_{0}$ as
$k\to\infty$, and $\{M_{k}\}$ satisfies
$\lim_{k\to\infty}k^{-N/2}M_{k}=\infty$. Furthermore it can be
assumed that $\norm{\eta_{k}}_{L^2}\leq c_{0}k^{N/2}$. The
single-point initial blow-up will be a consequence of
 \blemma {Energ}
For any $\delta>0$ there exists $C=C(\delta)$ such that:
  \begin{equation}\label{energ2}
  \sup_{t\in[0,1]}\myint{\abs x\geq\gd}{}u_{k}^2(x,t)dx
  +\myint{0}{1}\myint{\abs x\geq\gd}{}(\abs{\nabla
u_k}^2+u_{k}^2)dx\,dt\leq C\forevery k\in\BBN.
 \end {equation}
 \es
\Proof For $r\in (0,1)$, $\gt\geq 0$ we set $\Gw(\gt)=\{x\in\BBR^N:\abs x>\gt\}$,
$Q^r(\gt)=\Gw(\gt)\ti (0,r]
$, $Q_r(\gt)=\Gw(\gt)\ti(r,1)$  and $Q_r=\BBR^N\ti(r,1)$, and
denote
 $$I_1(r)=\myint{}{}\!\!\!\myint{Q_r}{}\abs {\nabla u}^2dx\,dt,\;
 I_2(r)=\myint{}{}\!\!\!\myint{Q_r}{}u^2dx\,dt,
 I_3(r)=\myint{}{}\!\!\!\myint{Q_r}{}h(t)\abs u^{q+1}dx\,dt.
 $$
 If we multiply the equation by $u(x,t)e^{(r-t)/(2-r)}$, integrate on $Q_r$ and use H\"older's
 inequality, we get, since $h$ is nondecreasing,
   \begin{multline}\label{energ4}
\myint{\BBR^N}{}u^2(x,1)dx+I_1(r)+I_2(r)+I_3(r)\leq c\int_{\mathbb{R}^N}u^2(x,r)\,dx\\
\leq c\gt^{\frac{N(q-1)}{q+1}}h(r)^\frac{-2}{q+1}
(-I'_3(r))^\frac{2}{q+1}+c\myint{\Gw(\gt)}{}u^2(x,r)dx.
 \end{multline}
Let $\gt\mapsto\gm(\gt)$ be a smooth decreasing function, we define
  $$\BA {l}
   E_{1}^\gm(r,\gt)=\myint{}{}\!\!\!\myint{Q^r(\gt)}{}\!\!\!\left(\abs
{\nabla u}^2
 +\gm^2u^2(x,t)\right)e^{-\gm^2t}dx\,dt,\\[4mm]
E_2(r,\gt)=\myint{}{}\!\!\!\myint{Q^r(\gt)}{}\!\!\!u^2dx\,dt\quad\mbox
  {and }\;
 f_\gm(r,\gt)=\sup\{e^{-\gm^2t}\myint{\Gw(\gt)}{}\!\!\!u^2(x,t)dx:0\leq t\leq r\}
 \EA$$

\noindent and $f(r)= f_0(r,0)$.  Then we introduce a parameter in the equation as in  \cite {OR} by multiplying it by
$u(x,t)\exp(-\mu^2(\tau)t)$ and integrating in the domain $Q^r(\tau)$ with $\tau>k^{-1}$
$Q^r(\tau)$ and $\tau>k^{-1}$. After some simple computations we deduce:
 $$  f_\gm(r,\gt)+ 2 E_{1}^\gm(r,\gt)
 \leq \myfrac {2}{\gm}
 \myint{0}{r}\myint{\prt \Gw(\gt)}{}\!\!\!\left(\abs {\nabla u}^2
 +\gm^2u^2(x,t)\right)e^{-\gm^2t}dS\,dt\qquad\forall\,\tau>k^{-1}.
 $$
Assuming $1-2\gm'/\gm^2>1/2$, we deduce from last inequality:
 $$ f_\gm(r,\gt)+ E_{1}^\gm(r,\gt)\leq
 -\myfrac {2}{\gm(\tau)}\myfrac {d E_{1}^\gm(r,\gt)}{d\gt} \qquad\forall\,\tau>k^{-1},
 $$
 and by integration
$$f_\gm(r,\gt_2)\left(e^{\int_{\gt_1}^{\gt_2}\frac{\gm(\gt)\,d\gt}{2}}-1\right)
+ E_{1}^\gm(r,\gt_{2})e^{\int_{\gt_1}^{\gt_2}\frac{\gm(\gt)\,d\gt}{2}}
 \leq  E_{1}^\gm(r,\gt_{1})\forevery\gt_2>\gt_1>k^{-1}.
 $$
The choice $\gm(\gt)=r^{-1}(\gt-k^{-1})/8$ ($\gt>k^{-1}$) yields to
\begin{multline}\label{16}
\myint{\Gw(\gt)}{}\!\!\!u^2(x,r)dx +
\myint{}{}\!\!\!\myint{Q^r(\gt)}{}\!\!\!\left(\abs{\nabla_x u}^2+
\frac {(\gt-k^{-1})^2}{64r^2}u^2\right)dx\,dt \leq c_1e^{-\frac
{(\gt-k^{-1})^2}{64r}}\\\times
\myint{}{}\!\!\!\myint{Q^r(\gt_0^k)}{}\!\!\!\left(\abs{\nabla_x
u}^2+ \frac {u^2}{2r}\right)dx\,dt\forevery \gt\geq
\tilde\gt_0^k:=k^{-1}+8\sqrt{r}>\tau_0^k:=k^{-1}+4\sqrt{2r}.
\end{multline}

\noindent We will need standard global energy estimate of solution
of problem \eqref{energ1} too:
\begin{multline}\label{17}
\int_{\mathbb{R}^N}|u(x,r)|^2\,dx+\int_{Q^r}(|\nabla_xu|^2+|u|^2+h(t)|u|^{q+1})\,dx\,dt\\
\leq c\|u_{0,k}\|^2_{L_2(\mathbb{R}^N)}\leq \bar
cM_k\qquad\forall\,r>0.
\end{multline}
Estimating the right-hand side terms in \eqref{energ4} and \eqref{16} by
\eqref{17}, we derive:
 \begin{equation}\label {energ5}\BA {l}
(i)\;\displaystyle\sum_{i=1}^{3}I_{i}(r) \leq
c_{1}\gt^{\frac{N(q-1)}{q+1}}h(r)^\frac{-2}{q+1}
(-I'_3(r))^\frac{2}{q+1}+\myfrac {c_{2}M_{k}}{r}e^{-(\gt
-k^{-1})^{2}/64r}\quad\forall\, \tau\geq\tilde\tau_0^k(r)\\
(ii)\; f_{0}(r,\gt)+E^{0}_{1}(r,\gt)+\myfrac {\gt -k^{-1}}{64r^{2}}E_{2}(r,\gt)
 \leq \myfrac {c_{2}M_{k}}{r}e^{-(\gt -k^{-1})^{2}/64r}\qquad\forall\, \tau\geq\tilde\tau_0^k(r).\!\!\!
 \EA\end{equation}
  Next we choose $M_{k}=e^{e^k}$, fix $\ge_{0}\in (0,e^{-1})$ and define a pair $(r_{k},\ \gt_k)$
  by the following relations:
 $r_{k}=\sup\{r:I_{1}(r)+I_{2}(r)+I_{3}(r)>2 M^{\ge_{0}}_{k}\}$; $c_2r_k^{-1}\exp(-\frac{\gt_k^2}{64r_k})M_k=M_k^{\ge_0}$
 $\Leftrightarrow$ $\gt_{k}=8\sqrt{r_{k}(1-\ge_{0})e^k+\ln (c_{2}/r_{k})}$.
 Taking $\gt=\gt_{k}+k^{-1}$ in (\ref {energ5})-i and solving the corresponding O.D.E. yields
 the estimate:
  \begin{equation}\label {energ6}\BA {l}
\displaystyle\sum_{i=1}^{3}I_{i}(r)\leq
c_{3}(\gt_{k}+k^{-1})(H(r))^{-2/(q-1)}\forevery r\leq r_{k} \mbox
{ , }H(r)=\myint{0}{r}h(s)ds. \EA\end{equation} If we write
$h(t)=e^{-\gw(t)/t}$, the assumption $\inf\{\gw(t)/t^\ga:0<t\leq
1\}>0$ implies that $H(r)\geq  c_0e^{-\gw(r)/r}r^2/\gw(r)$ and,
replacing $\gt_{k}$ by its expression, (\ref {energ6}) turns into
$$\BA {l}
\displaystyle\sum_{i=1}^{3}I_{i}(r)\leq
c_{4}(\sqrt{r_{k}(1-\ge_{0})e^k+\ln
(c_{2}/r_{k})}+k^{-1})^N\left(\myfrac{\gw(r)e^{-\gw(r)/r}}{r^{2}}\right)^{2/(q-1)}\!\!\!\!
\forall r\leq r_{k}. \EA$$ Thus $r_{k}\leq b_{k}$, where $b_k$ is
solution of equation:
$$c_{4}\left(\sqrt{r_{k}(1-\ge_{0})e^k+\ln (c_{2}/b_{k})\;}+k^{-1}\right)^N\!\!\!\left(\myfrac{\gw(b_{k})e^{-\gw(b_{k})/b_{k}}}{b_{k}^{2}}\right)^{2/(q-1)}\!\!\!\!\!\!=2M_{k}
^{\ge_{0}}=2e^{\ge_{0}e^k}.
$$
 From this inequality using additionally assumption on $\gw(t)$,
 we obtain inequalities:
$c_{5}e^k\geq \gw(b_{k})/b_{k}\geq c_{6}e^k,\ c_6>0$; $b_{k}\geq
e^{-c_{7}k}, \ c_7>0$. These inequalities yield:
  \begin{equation}\label {energ7}\BA {l}
\gt_{k}\leq c_{8}\sqrt{\gw(c_{9}e^{-k})}\;. \EA
\end{equation}
Using the definition of $r_{k}$, inequality (\ref{energ5})-ii, the
fact that $3M_{k}^{\ge_{0}}\leq\bar c M_{k-1}\ \forall\,k\geq
k_0(\bar c)$ ($\bar c$ is from \eqref{17}, $0<\ge_0<e^{-1}$), we
deduce the main result of first round of computations:
  \begin{equation}\label {energ8}\BA {l}
\displaystyle\sum_{i=1}^{3}I_{i}(r_{k})+f_{0}(r_{k},\gt_{k}+k^{-1})
+\sum_{i=1}^{2}E_{i}(r_{k},\gt_{k}+k^{-1})\leq
3M_{k}^{\ge_{0}}\leq \bar cM_{{k-1}}. \EA
\end{equation} Next we
organize the second round of estimates with
$\gm(\gt)=(\gt-\gt_{k}-k^{-1})/8$, $r_{k-1}$ and $\gt_{k-1}$ be
defined similarly as $r_{k}$ and $\gt_{k}$, up to the change of
indices, using obtained estimate \eqref{energ8} instead of
\eqref{17}. As result we derive:
  \begin{equation}\label {energ9}
\sum_{i=1}^{3}I_{i}(r_{k-1})\!+\!f_{0}(r_{k-1},\gt_{k}\!+\gt_{k-1}\!+\!k^{-1})
\\\!+\sum_{i=1}^{2}E_{i}(r_{k-1},\gt_{k}\!+\gt_{k-1}\!+\!k^{-1})\leq \bar
cM_{{k-2}}.
\end{equation}
Fixing arbitrary $n>k_0(\bar c)$ and repeating the above described round of
computations $k-n$ times, we obtain:
  \begin{equation}\label {energ9}\BA {l}
\displaystyle\sum_{i=1}^{3}I_{i}(r_{n})+f_{0}\left(r_{n},\sum_{j=0}^{k-n}\gt_{k-j}+k^{-1}\right)
+\sum_{i=1}^{2}E_{i}\left(r_{n},\sum_{j=0}^{k-n}\gt_{k-j}+k^{-1}\right)\leq
\bar cM_{n-1}, \EA
\end{equation} and, since by induction $\gt_{k-j}$
satisfies (\ref{energ7}) with $k$ replaced by $k-j$, we obtain
  \begin{equation}\label {energ10}\BA {l}
\displaystyle\sum_{j=0}^{k-n}\gt_{k-j}\leq
c_{8}\sum_{j=0}^{k-n}\sqrt{\gw(c_{9}e^{-(k-j)})} \leq
c_{10}\myint{c_{9}e^{-k}}{c_{9}e^{-n}}\myfrac{\sqrt{\gw(s)}ds}{s}.
\EA\end{equation} We denote
$\gt^*(n)=\lim_{k\to\infty}c_{8}\sum_{j=0}^{k-n}\sqrt{\gw(c_{9}e^{-(k-j)})}$. We
derive from (\ref{energ9}) by letting $k\to\infty$,
  \begin {equation}\label {energ11}\BA {l}
\displaystyle\sup_{0<t\leq r_{n}}\myint{\abs
x\geq\gt^*(n)}{}\!\!\!u^2(x,t)dx+\myint{0}{r_{n}}\myint{\abs
x\geq\gt^*(n)}{}\!\!\!(\abs{Du}^2+u^2)dx\,dt \leq \bar cM_{n-1}.
\EA\end {equation} 
Due to assumption \eqref{main2} $\tau^*(n)\to
0$ as $n\to\infty$, therefore inequality \eqref{energ11} implies
the result.


\end {document}